\documentclass[11pt]{amsart}

\usepackage{amssymb,amsmath,amsthm,amscd}

\input{amssym.def}
\newcommand{\bt}{\begin{Theorem}}
\newcommand{\et}{\end{Theorem}}
\newcommand{\bi}{\begin{itemize}}
\newcommand{\ei}{\end{itemize}}
\newcommand{\bea}{\begin{eqnarray}}
\newcommand{\eea}{\end{eqnarray}}

\theoremstyle{plain}

\newtheorem{Theorem}{\sc Theorem}[section]
\newtheorem{Lemma}{\sc Lemma}[section]
\newtheorem{Proposition}{\sc Proposition}[section]

\theoremstyle{definition}
\newtheorem{Definition}{\sc Definition}[section]

\theoremstyle{remark}
\newtheorem{Remark}{\sc Remark}[section]

\newcommand{\be}{\begin{equation}}
\newcommand{\ee}{\end{equation}}

\newcommand{\newsection}{\setcounter{equation}{0}\setcounter{Theorem}{0}
\setcounter{Proposition}{0}\setcounter{Lemma}{0}\setcounter{Corollary}{0}
\setcounter{Definition}{0}\setcounter{Example}{0}\setcounter{Question}{0}
\setcounter{Remark}{0}}

\newcommand{\chd}{Cowen - Douglas }%
\def\hom{{homogeneous~}}%
\def\op{{operator~}}%
\def\Op{{operator}}%
\def\ops{{operators~}}%
\def\Ops{{operators}}%
\def\rep{{representation~}}%
\def\reps{{representations~}}%
\def\Rep{{representation}}%
\def\mult{{multiplier~}}%
\def\irr{{irreducible~}}%
\def\hs{{Hilbert space~}}%
\newcommand{\rk}{reproducing kernel }%

\newcommand{\bm}[1]{\mathbf #1}
\newcommand{\bb}[1]{\mathbb #1}
\newcommand{\cl}[1]{\mathcal #1}
\newcommand{\ger}[1]{\mathfrak #1}
\newcommand{\scr}[1]{\mathcal #1}

\def\C{\mbox{${\mathbb C}$}}
\def\D{\mbox{${\mathbb D}$}}

\def\N{\mbox{${\mathbb N}$}}

\def\R{\mbox{${\mathbb R}$}}

\newcommand{\g}{{\mathfrak g}}%

\newcommand{\m}{{\boldsymbol \mu}}%

\newcommand{\inner}[2]{\langle #1,#2 \rangle }%
\newcommand{\mat}[4]{\Big ( \begin{matrix} #1 & #2\\ #3 & #4 
\end{matrix} \Big )}%

\title[Homogeneous operators]{Homogeneous operators on Hilbert spaces of\\
holomorphic functions -- I}

\author{Adam Kor\'{a}nyi}
\address{Lehman College\\
The City University of New York\\
Bronx, NY 10468
}
\email{adam.koranyi@lehman.cuny.edu}
\author{Gadadhar Misra}
\address{Indian Statistical Institute\\
R.\ V.\ College Post\\
Bangalore 560 059
 }
\email{gm@isibang.ac.in}
\thanks{This work was supported in part by DST -
NSF S\&T Cooperation Programme}

\begin{document}
\baselineskip=14.25pt
\begin{abstract}

  In this paper we construct a large class of multiplication operators
  on \rk Hilbert spaces which are {\em homogeneous} with respect to
  the action of the M\"{o}bius group consisting of bi-holomorphic
  automorphisms of the unit disc $\D$.  For every $m\in \N$ we have a
  family of operators depending on $m+1$ positive real parameters.
  The kernel function is calculated explicitly.  It is proved that
  each of these operators is bounded, lies in the \chd class of $\bb
  D$ and is irreducible. These \ops are shown to be mutually pairwise
  unitarily inequivalent.
 
\end{abstract}
\maketitle
\section{Introduction} 
A \hom \op on a \hs $\scr H$ is a bounded \op $T$ whose spectrum is
contained in the closure of the unit disc $\bb D$ in $\bb C$ and is
such that $g(T)$ is unitarily equivalent to $T$ for all linear
fractional transformations $g$ which map $\bb D$ to $\bb D$.  This
class of \ops has been studied in a number of articles \cite{surv,
  CM, con, W, shift, Peng-Z, Arazy-Z, QME2}. It is known that every
\hom \op is a block shift, that is, $\scr H$ is the orthogonal direct
sum of subspaces $V_n$, indexed by all integers, all non-negative
integers or all non-positive integers, such that $T(V_n) \subseteq
V_{n+1}$ for each $n$. 

The case where $\dim V_n = 1$ for each $n$ is completely known, the
corresponding \ops have been classified in \cite{shift}.  The
classification in the case where $\dim V_n \leq 2$ and $T$ belongs to
the \chd class of $\bb D$ is complete and the \ops are explicitly
described in \cite{W}.  Beyond this there are only some results of a
general nature, and not too many examples are known (cf. \cite{surv}).

In the present article we construct a large family of examples. For every 
natural number $m$ we construct a family depending on $m+1$ parameters. 
Each one of the examples is realized as the multiplication \op on a \rk 
space of vector-valued holomorphic functions.
All of these \rk Hilbert spaces admit a direct sum decomposition 
$\oplus_{n\geq 0} V_n$ with $\dim V_n =  n+1$ if $0 \leq n < m$ and 
$\dim V_n = m+1$ for $n \geq m$.  The reproducing kernels are described 
explicitly.  All our examples are \irr \ops and their adjoints belong to the 
\chd class.

We have chosen a presentation as elementary as possible, based on explicit 
computations. This seemed to be appropriate here since our goal was a 
complete explicit description of the examples.  On the other hand, it
does not explain the deeper background of the results.  To remedy this 
situation we have added a final section which discusses a more conceptual 
approach to the examples.  In a planned expository article on the subject 
there will be more details about the various ways in which one can arrive at 
the construction of our examples.    
    
The more conceptual approach will play a leading role in the sequel to the 
present article, where a description of all \hom \chd \ops will be given 
albeit in a less explicit way than our present examples.

Our results are also the subject of a short note presented to the 
Comptes Rendus de l'Acad\'{e}mie des Sciences, Paris \cite{KM}. 
\newsection
\section{Preliminaries}
We denote by $\bb D$ the open unit disc in $\bb C$ and by $G$ the
group of M\"{o}bius transformations $z \mapsto \tfrac{a z +
b}{\bar{b}z+\bar{a}}$, $|a|^2 - |b|^2 = 1$.  
Let $G_0$ be the group ${\rm SU}(1,1) = 
\Big \{ \Big (\begin{matrix} a & b\\ \bar{b}& \bar{a}
\end{matrix} \Big ): |a|^2 - |b|^2 = 1  \Big \}$.
So, $G= G_0/\{\pm I\}$.  We 
denote by $\tilde{G}$, the universal covering group of $G$. 


All Hilbert spaces $\scr H$ considered in this article will be spaces 
of holomorphic functions $f:\bb D \to V$ taking their values in a
finite dimensional Hilbert space $V$ and possessing a reproducing
kernel $K$.  A \rk is a function $K: \bb D \times \bb D \to {\rm
Hom}(V,V)$ holomorphic in the first variable and anti-holomorphic in
the second, such that $K_\omega \zeta$ defined by $(K_\omega\zeta)(z) 
:= K(z,\omega) \zeta$
is in $\scr H$ for each $\omega\in \bb D$, $\zeta\in V$, and 
\begin{equation}  
\inner{f}{K_\omega\zeta}_{\mathcal H} = \inner{f(\omega)}{\zeta}_V
\end{equation}
for all $f\in \scr H$.

As is well known, if $\{e_n\}_{n=0}^\infty$ is any orthonormal basis of 
$\scr H$, then we have 
\begin{equation} \label{Kexpand}
K(z,\omega) = \sum_{n=0}^\infty e_n(z) e_n(\omega)^*
\end{equation}
with the sum converging pointwise.  Here we interpret a formal product
$\xi \eta^*$ for $\xi, \eta \in V$ as the transformation $\zeta\mapsto
\inner{\zeta}{\eta}\xi$; when $V=\bb C^k$, $k \in \bb N$, and its elements 
are written as column vectors, $\xi \eta^*$ is just the usual matrix product.

We will be concerned with \mult \reps of $\tilde{G}$ on the Hilbert
space $\scr H$.  A {\em multiplier} is a continuous function $J:\tilde{G}
\times \bb D \to {\rm Hom}(V,V)$, holomorphic on $\mathbb D$, such that 
\begin{equation} \label{mult}
J(gh, z) = J(h,z)J(g, hz)
\end{equation}
for all $g,h \in \tilde{G}$ and $z \in \bb D$.  For $g\in\tilde{G}$, 
we define $U(g)$ on  on ${\rm Hol}(\bb D,V)$ by 
\begin{equation} \label{Udef}
(U(g) f)(z) = J(g^{-1}, z) f (g^{-1}(z)). 
\end{equation}
It is easy to see that the multiplier identity (\ref{mult}) is equivalent 
to $U(gh) = U(g)U(h)$. 

Suppose that the action $g\mapsto U(g)$, $g\in \tilde{G}$,
defined in (\ref{Udef}) preserves $\scr H$ and is unitary on it, then
we say that $U$ is a {\em unitary multiplier representation} of
$\tilde{G}$.

Also, if the reproducing kernel $K$ transforms according to the rule 
\begin{equation} \label{Ktransrule}
J(g,z)K(g(z), g(\omega)) J(g,\omega)^* = K(z,\omega)
\end{equation}
for all $g \in \tilde{G}$; $z,\omega \in \bb D$,  then we say that $K$ is 
{\em quasi-invariant}.

\begin{Proposition} \label{}
Suppose $\scr H$ has a reproducing kernel $K$.  Then $U$ defined by 
(\ref{Udef}) is a unitary representation if and only if $K$ is quasi-invariant.
\end{Proposition}
\begin{proof}
Assume that $K$  is quasi-invariant. We have to show that the linear 
transformation $U$ defined in (\ref{Udef}) is unitary. We note, 
writing $\tilde{\omega} = g^{-1}(\omega)$ and $\tilde{\omega}^\prime 
= g^{-1}(\omega^\prime)$,  
\begin{align*}
\inner{U(g^{-1})K(\cdot, \omega)\xi}{U(g^{-1})K(\cdot, \omega^\prime)\eta} &=
\inner{J(g,\cdot) K(g(\cdot), \omega)\xi}{J(g,\cdot)
K(g(\cdot), \omega^\prime)\eta}\\   
& = \inner{K(\cdot,\tilde{\omega})J(g,\tilde{\omega})^{*-1}\xi}
{K( \cdot,\tilde{\omega}^\prime) J(g,\tilde{\omega}^\prime)^{*-1}\eta}\\
&=\inner{K(\tilde{\omega}^\prime,\tilde{\omega})J(g,\tilde{\omega})^{*-1}\xi}
{J(g,\tilde{\omega}^\prime)^{*-1}\eta}\\
&= \inner{J(g,\tilde{\omega}^\prime)^{-1}K(\tilde{\omega}^\prime,
\tilde{\omega})J(g,\tilde{\omega})^{*-1}\xi}{\eta}\\
&= \inner{K(\omega^\prime, \omega)\xi}{\eta} 
\end{align*}
and it follows that $U(g^{-1})$ is isometric.  

On the other hand, if $U$ of (\ref{Udef}) is unitary then 
the reproducing kernel $K$ of the  \hs $\mathcal{H}$ satisfies
the transformation rule (\ref{Ktransrule}).
A reproducing kernel $K$ has the expansion (\ref{Kexpand}).
It follows from the uniqueness of the \rk that the expansion is
independent of the choice of the orthonormal basis.  Consequently, we
also have $K(z,\omega) = \sum_{\ell=0} (U_{g^{-1}}e_\ell)
(z)(U_{g^{-1}}e_\ell) (\omega)^*$ which verifies the equation 
(\ref{Ktransrule}). 
\end{proof}

When we are in the situation of the Proposition and if we can prove that
the operator $M$ defined by $(Mf)(z)=z f(z)$ is 
bounded on $\scr H$, then $M$ is a homogeneous operator. 
This is well-known and trivial: 
Clearly, $(g(M)f)(z) = g(z)f(z)$ and hence   
$(M U(g^{-1})f )(z)= z J(g,z) f(g(z)) = J(g,z) g^{-1}(g(z)) f(g(z)) 
= (U(g^{-1}) (g^{-1}(M)) f)(z)$, for all $g\in \tilde{G}$, 
$f\in \scr H$, $z\in \bb D$.  If, in addition, $\dim \ker (M-\omega I)^* =n$ 
and the operator $(M-\omega I)^*$ is bounded below, on the orthogonal complement of its kernel, 
for every $\omega \in \bb D$ then $M^*$ is in the Cowen-Douglas class 
(see \cite{C-D}) $B_n(\bb D)$. 

In the case of \rk Hilbert spaces of scalar functions (i.e. when $\dim
V=1$) the unitary multiplier \reps of $\tilde{G}$ are well-known.  We
describe them here because they will be used in the next section.
They are the elements of the holomorphic discrete series depending on
one real parameter $\lambda > 0$. They act on the Hilbert space
$A^{(\lambda)}(\bb D)$ characterized by its \rk $B^\lambda(z,\omega)=
(1-z\bar{\omega})^{-2\lambda}$. Here
$B(z,\omega)=(1-z\bar{\omega})^{-2}$ is the \rk of the Bergman space
$A^2(\bb D)$, the Hilbert space of square integrable 
(with respect to normalized area measure) holomorphic functions on the unit
disc $\bb D$.

For $g\in \tilde{G}$, $g^\prime(z)^\lambda$ is a real analytic function on the 
simply connected set $\tilde{G}\times \D$, holomorphic in $z$. 
Also $g^\prime(z)^\lambda \not = 0$ since $g$ is one-one and holomorphic.  
Given any $\lambda \in \C$, taking the principal branch of the power function 
when $g$ is near the identity, we can uniquely define $g^\prime(z)^\lambda$ 
as a real analytic function on $\tilde{G} \times \D$ which is holomorphic on 
$\D$ for all fixed $g\in \tilde{G}$.  The multiplier
$j_\lambda(g,z)=g^\prime(z)^\lambda$ defines on
$A^{(\lambda)}(\bb D)$ the unitary \rep $D_\lambda^+$ by the formula 
(\ref{Udef}), that is, 
\begin{equation} \label{defn disc repn}
D_{\lambda}^+ (g^{-1})(f) =
(g^\prime)^{\lambda} (f \circ g),\: f\in A^{(\lambda)}(\bb D),\:g\in \tilde{G}.
\end{equation}
An orthonormal basis of the space is given by
$\Big\{\sqrt{\tfrac{(2\lambda)_n}{n!}}z^n\Big \}_{n \geq 0}$, where 
$(x)_n=x(x+1)\ldots (x+n-1)$ is the 
Pochhammer symbol. The operator $M$ is bounded on the Hilbert space 
$A^{(\lambda)}(\bb D)$.  
It is easily seen to be in the Cowen-Douglas class $B_1(\bb D)$.
\newsection
\section{Construction of the Hilbert spaces and representations}

Let $\mathrm{Hol}(\bb D, \C^k)$ denote the vector space of
all holomorphic functions on $\D$ taking values in $\C^k$, $k\in \bb N$.  
Let $\lambda$ be a real number and $m$
be a positive integer satisfying $2\lambda - m > 0$. For brevity, we will     
write $2\lambda_j = 2\lambda-m+2j$. 

For each $j$, $0 \leq j \leq m$, define the operator 
$\Gamma_j:A^{(\lambda_j)}(\D) \to \mathrm{Hol}(\bb D, \C^{m+1})$ by the formula 
$$
(\!\!(\Gamma_jf)\!\!)_\ell = \begin{cases}
\tbinom{\ell}{j}\frac{1}{(2\lambda_j)_{\ell-j}} f^{(\ell -j)} &
\mbox{~if~} \ell \geq j\\ 0 &\mbox{~if~} \ell < j, 
\end{cases} 
$$
for $f\in A^{(\lambda_j)}(\D), \, 0\leq \ell \leq m.$ 
Here $(\!\!( \Gamma_jf )\!\!)_\ell$ denotes the $\ell$th component of the 
function $\Gamma_jf$ and $f^{(\ell -j)}$ denotes the $(\ell-j)$th 
derivative of the holomorphic function $f$.

We denote the image of $\Gamma_j$ by $\bm{A}^{(\lambda_j)}(\D)$ and
transfer to it the inner product of $A^{(\lambda_j)}(\D)$, that is, we
set $\inner{\Gamma_j f}{\Gamma_j g} = \inner{f}{g}$, for $f,g\in
A^{(\lambda_j)}(\D)$.  The \hs $\bm{A}^{(\lambda_j)}(\D)$ is a \rk
space because the point evaluations $f \mapsto (\Gamma_jf)(\omega)$
are continuous for each $\omega \in \bb D$.  Let
$\bm{B}^{(\lambda_j)}$ denote the reproducing kernel for the \hs
$\bm{A}^{(\lambda_j)}(\D)$.

The algebraic sum of the linear spaces $\bm{A}^{(\lambda_j)}(\D)$,
$0\leq j \leq m$ is direct. This is easily seen.  If $\sum_{j=0}^m
\Gamma_j f_j = 0$, $f_j \in A^{(\lambda_j)}(\D)$, then $f_0 =
(\!\!(\Gamma_0f_0 )\!\!)_0 = 0$ since $(\!\!(\Gamma_jf_j )\!\!)_0 = 0$ for $j >0$.
Similarly, $f_1 = (\!\!(\Gamma_1f_1)\!\!)_1 = 0$ since $(\!\!(\Gamma_jf_j)\!\!)_1 = 0$
for $j > 1$.  Continuing in thpositive numbersis fashion, we see that $f_m =0$. It
follows that we can choose $m$ positive numbers, 
$\mu_j,\,1\leq j \leq m$, set $\mu_0=1$, write   
$\m= (\mu_0, \mu_1, \ldots , \mu_m)$, and define an inner
product on the direct sum of the $\bm A^{(\lambda_j)}(\bb D)$ by
setting
\begin{equation}
\inner{\sum_{j=o}^m \Gamma_jf_j}{ \sum_{j=o}^m \Gamma_j g_j} = 
\sum_{j=0}^m \mu_j^2 \inner{f_j}{g_j}, 
\:\:f_j,g_j\in A^{(\lambda_j)}.  
\end{equation}
We obtain a Hilbert space in this manner which we denote by 
$\bm{A}^{(\lambda, \m)}(\D)$.
It has the \rk $\bm{B}^{(\lambda,\m)} = \sum_{j=0}^m \mu_j^2 \,
\bm{B}^{(\lambda_j)}$.

The direct sum of the discrete series \reps $D_{\lambda_j}^+$ on
$\oplus_{j=0}^m A^{(\lambda_j)}$ can be transferred to $\bm
A^{(\lambda, \m)}(\D)$ by the map $\Gamma = \oplus_{j=0}^m \mu_j
\Gamma_j$.  It is a unitary \rep of the group $\tilde{G}$ which we
call $U$.  Its \irr subspaces are the $\bm A^{(\lambda_j)}(\bb D)$.

We will show that $U$ is a \mult \Rep.  For each $\bm
A^{(\lambda_j)}(\bb D)$ separately this is fairly obvious by checking
the effect of $\Gamma_j$.  The important point is that the multiplier is
the same on each $\bm A^{(\lambda_j)}(\bb D)$.

We need a relation between $g^{\prime \prime}(z)$ and $g^\prime(z)$. The elements of 
$G_0$ are the matrices $\mat{a}{b}{\bar{b}}{\bar{a}}$, $|a|^2 - |b|^2 =1$, acting on 
$\D$ by fractional linear transformations.  The inequalities 
\begin{equation} \label{ineq}
|a-1| < 1/2, \,\, |b| < 1/2
\end{equation}
determine a simply connected neighborhood $U_0$ of $e$ in $G_0$.  
Under the natural projections, it is diffeomorphic with a neighborhood $U$ 
of $e$ in $G$ and with a neighborhood $\tilde{U}$ of $e$ in $\tilde{G}$. 
So, we may use $a,\, b$ satisfying \eqref{ineq} to parametrize $\tilde{U}$.
For $g\in \tilde{U}$, $z\in \D$ we have $g^\prime(z)=(\bar{b} z + \bar{a})^{-2}$ and  
$g^{\prime \prime}(z)= -2 \bar{b}(\bar{b} z + \bar{a})^{-3}$, which gives a relation
\begin{equation}\label{doubleder}
g^{\prime \prime}(z)=-2cg^\prime(z)^{3/2},
\end{equation} 
where $c=c_g$ depends on $g$ real analytically and is independent of $z$; the 
meaning of $g^\prime(z)^{3/2}$ is as defined earlier. Since both sides are real 
analytic, \eqref{doubleder} remains true on all of $\tilde{G} \times \D$.

\begin{Definition}
Let $J: \tilde{G} \times \D \to \C^{m+1\times m+1}$ be the function 
given by the formula 
\begin{equation} \label{Jmult}
J(g, z)_{p,\ell} =\begin{cases} 
\tbinom{p}{\ell} (-c)^{p-\ell}
(g^\prime)^{\lambda-\frac{m}{2} + \frac{p + \ell}{2}}(z)&\mbox{~if~} p \geq \ell\\ 
0&\mbox{~if~} p < \ell, \end{cases}
\end{equation}
for  $g \in \tilde{G}$. Here $c$ is the constant depending on $g$ 
as in (\ref{doubleder})
\end{Definition}

The following Lemma is used for showing that  $U$ 
is a multiplier representation.
\begin{Lemma}
For any $g\in \tilde{G}$, we have the formula  
$$
\big ( (g^\prime)^\ell(f\circ g) \big )^{(k)} = 
\sum_{i=0}^k
\tbinom{k}{i}(2\ell+i)_{k-i}(-c)^{k-i}(g^\prime)^{\ell+\frac{k+i}{2}}\big ( 
f^{(i)} \circ g \big ). 
$$ 
\end{Lemma}

\begin{proof}
The proof is by induction, using the formula (\ref{doubleder}). 
For $k=0$, the formula is an identity.  Assume the formula to be valid for 
some $k$. Then 
\begin{eqnarray*}
\lefteqn{\big ( (g^\prime)^\ell(f\circ g) \big )^{(k+1)}}\\ 
& = &\sum_{i=0}^k \tbinom{k}{i}(2\ell+i)_{k-i}(-c)^{k-i} \Big 
\{(\ell+\frac{k+i}{2}) (g^\prime)^{\ell+\frac{k+i}{2} - 1} 
g^{\prime\prime} \big ( f^{(i)} \circ g \big )
+(g^\prime)^{\ell+\frac{k+i}{2}}\big (f^{(i+1)}\circ g \big ) 
g^\prime \Big \} \\
& = &\sum_{i=0}^k \tbinom{k}{i}(2\ell+i)_{k-i}(-c)^{k-i} \Big \{
(2\ell+k+i) (-c) (g^\prime)^{\ell+\frac{k+i+1}{2}}
\big ( f^{(i)} \circ g \big ) + (g^\prime)^{\ell+\frac{k+i+2}{2}}
\big (f^{(i+1)}\circ g \big ) \Big \}\\
&=&\sum_{i=0}^k \tbinom{k}{i}(2\ell+i)_{k-i}(2\ell+k+i)(-c)^{k+1-i}
(g^\prime)^{\ell+\frac{k+i+1}{2}}\big ( f^{(i)} \circ g \big )\\
& & \quad \quad \quad \quad  \quad \quad \quad \quad 
+ \sum_{i=1}^{k+1} \tbinom{k}{i-1}(2\ell+i-1)_{k+1-i}
(-c)^{k+1-i}(g^\prime)^{\ell+\frac{k+i+1}{2}}
\big ( f^{(i)} \circ g \big ).
\end{eqnarray*}
Now, we observe that 
\begin{eqnarray*}
\lefteqn{\tbinom{k}{i}(2\ell+i)_{k-i}(2\ell+k+i) + 
\tbinom{k}{i-1}(2\ell+i-1)_{k+1-i}} \quad \quad \quad\quad \quad \quad\\
& = & (2\ell+i)_{k-i} \big \{ \tbinom{k}{i}(2\ell+k+i) + 
\tbinom{k}{i-1} (2\ell + i -1) \big \} \\
& = &(2\ell+i)_{k-i}\big \{\big ( \tbinom{k}{i} + \tbinom{k}{i-1} \big )
(2\ell+k) +  i  \tbinom{k}{i} + (i-1+k) \tbinom{k}{i-1} \big \} \\
& = &(2\ell+i)_{k+1-i} \tbinom{k+1}{i}.
\end{eqnarray*}
Thus $\big ( (g^\prime)^\ell(f\circ g) \big )^{(k+1)}=
(2\ell+i)_{k+1-i} \tbinom{k+1}{i}(-c)^{k+1-i}
(g^\prime)^{\ell+\frac{k+i+1}{2}}$ completing the induction step.
\end{proof}
We can now prove the main theorem 
of this section.
\begin{Theorem} \label{mainthmsec3}
The image of $\oplus_{0}^m D_{\lambda_j}^+$ under $\Gamma$ is a
multiplier \rep with the multiplier given by $J(g, z)$ as in (\ref{Jmult}).
\end{Theorem}
\begin{proof} It will be enough to show 
$$
\Gamma_j \big ( D_{\lambda_j}^+(g^{-1}) f \big ) = 
J(g, \cdot)\big( (\Gamma_jf)\circ g \big )
$$
for each $j$, $0\leq j \leq m$.  We compute the $p$'th component 
on both sides.  

For $p < j$, both sides are zero by definition of $\Gamma_j$ and 
knowing that $J(g, z)_{p,\ell}=0$ for $\ell > p$.  For $p \geq j$, 
we have using the Lemma,  
\begin{align*}
(\!\!( (\Gamma_jD_{\lambda_j}^+(\varphi^{-1}) f)\big )\!\!)_p
&= \tbinom{p}{j}\frac{1}{(2\lambda_j)_{p-j}} 
\big ( (g^\prime)^{\lambda_j} f\circ g \big )^{p-j}\\
&=  \tbinom{p}{j}\frac{1}{(2\lambda_j)_{p-j}} \sum_{i=0}^{p-j} \tbinom{p-j}{i}
(2\lambda_j+i)_{p-j-i}(-c)^{p-j-i}(g^\prime)^{\lambda_j+\frac{p-j+i}{2}}
(f^{(i)}\circ g)\\
&= \tbinom{p}{j}\frac{1}{(2\lambda_j)_{p-j}} \sum_{\ell=j}^{p} 
\tbinom{p-j}{\ell -j}
(2\lambda_j+\ell -j)_{p-\ell}(-c)^{p-\ell}(g^\prime)^{\lambda_j - j +
\frac{p+\ell}{2}}(f^{(\ell-j)}\circ g)\\
&= \sum_{\ell=j}^m \tfrac{p!}{j!(\ell-j)!(p-\ell)!} 
\frac{1}{(2\lambda_j)_{\ell -j}}
(-c)^{p-\ell}(g^\prime)^{\lambda_j - j +\frac{p+\ell}{2}}
(f^{(\ell-j)}\circ g)\\
&=\sum_{\ell=0}^m J(\varphi,\cdot)_{p,\ell}
(\!\!( (\Gamma_jf)\circ g\big )\!\!)_\ell.
\end{align*}
\end{proof}

\newsection
\section{The orthonormal basis and the operator $M$}

The vectors $e_n^j(z) := \Gamma_j\big ( \sqrt{\tfrac{(2\lambda_j)_n}{n!}} 
z^n \big )$ clearly form an orthonormal basis in the Hilbert space 
$\bm A^{(\lambda_j)}(\D)$.  We have, by definition of $\Gamma_j$,   
\begin{equation} \label{onb bmA_j}
(\!\!( e_n^j(z))\!\!)_\ell = \begin{cases} 0 & \ell < j 
\mbox{~or~} \ell > n+j\\ 
\binom{\ell}{j}
\frac{\sqrt{n!}}{(n-\ell+j)!}\frac{\sqrt{(2\lambda_j)_n}}{(2\lambda_j)_{\ell-j}}
z^{n-\ell+j} & \ell \geq j \mbox{~and~}\ell \leq n+j.
\end{cases}
\end{equation}

We compute the reproducing kernel $\bm{B}^{(\lambda_j)}$ for the Hilbert space
$\bm{A}^{(\lambda_j)}(\D)$.  We have 
\begin{eqnarray} \label{Bergj}
\bm{B}^{(\lambda_j)}(z,\omega) &=& \sum_{n=0}^\infty
\big ( (\Gamma_j e^j_n)(z)\big ) \big ( (\Gamma_j e^j_n)(\omega)\big )^* 
\nonumber\\
&=&\big (\Gamma_j \sum_{n=0}^\infty e^j_n(z) \big ) \big (
\Gamma_j \sum_{n=0}^\infty {e^j_n(\omega)} \big )^*\nonumber\\
&=& \Gamma_j^{(z)} \Gamma_j^{(\bar{\omega})} B^{\lambda_j}(z,\omega), 
\end{eqnarray}
since the series converges uniformly on compact subsets. Explicitly,   
\begin{equation} \label{Bergman}
\bm B^{(\lambda_j)}(z,\omega)_{p,\ell} = \begin{cases}
\binom{\ell}{j} \binom{p}{j} \frac{1}
{(2\lambda_j)_{\ell-j}}\frac{1}{(2\lambda_j)_{p-j}} \partial^{(p-j)}
\bar{\partial}^{(\ell-j)}B^{\lambda_j}(z,\omega) & \mbox{if~} \ell,p \geq j\\
0 & \mbox{otherwise}.
\end{cases} 
\end{equation}
In particular, it follows that $\bm B^{(\lambda_j)}(0,0)$ is diagonal, and 
\begin{equation} \label{Bj(0,0)}
\bm B^{(\lambda_j)}(0,0)_{\ell,\ell} = \begin{cases} 0 & 
\mbox{if~} \ell < j\\
\binom{\ell}{j}^2 \frac{(\ell -j)!}{(2\lambda_j)_{\ell-j}} & \mbox{if~} 
\ell \geq j.
\end{cases}
\end{equation}
Then 
\begin{equation} \label{B(0,0)}
\bm B^{(\lambda,\m)}(0,0)_{\ell ,\ell} =  
\sum_{j=0}^m \bm B^{(\lambda_j)}(0,0)\,\mu_j^2. 
\end{equation}
A useful formula for $\bm B^{(\lambda, \m)}(z,\omega)$ can be easily obtained 
using \eqref{Ktransrule}.  
For $z\in \D$, we set $p_z=\frac{1}{\sqrt{(1-|z|^2)}}\mat{1}{z}{\bar{z}}{1}\in 
{\rm SU}(1,1)$. We also write $p_z$ for the corresponding element 
of $\tilde{G}$ such that $p_z$ depends contnuously on $z\in \D$ and $p_0 = e$. 
Then $p_z(0)=z$; $p_z^{-1}=p_{-z}$. By Theorem \ref{mainthmsec3}, formula \eqref{Ktransrule} 
holds for $\bm B^{\lambda, \m}$ and gives 
\begin{equation} \label{B(z,z)}
J_{p_{-z}}(z) \bm B^{\lambda, \m}(0,0) J_{p_{-z}}(z)^* = \bm B^{\lambda, \m}(z,z).
\end{equation}
We have 
$p_{-z}^\prime(\zeta)=\tfrac{1-|z|^2}{(1-\bar{z}\zeta)^2}$;  
$p_{-z}^\prime(z)=(1-|z|^2)^{-1}$. The  $-c$ of \eqref{doubleder} corresponding to 
$p_{-z}$ is $\tfrac{\bar{z}}{1-|z|^2}$.  So \eqref{Jmult} gives  
$$
J_{p_{-z}}(z)_{p,\ell} = \begin{cases}
(1-|z|^2)^{-\lambda-\tfrac{m}{2}}
\tbinom{p}{\ell} \bar{z}^{p-\ell}(1-|z|^2)^{m-p} & p \geq \ell \\
0 & p < \ell, 
\end{cases} 
$$
which can be written in matrix form as
\begin{equation} \label{Jexp}
J_{p_{-z}}(z) = (1-|z|^2)^{-\lambda-\tfrac{m}{2}}  D(|z|^2)\exp(\bar{z} S_m),  
\end{equation}
where $D(|z|^2)_{p,\ell} = (1-|z|^2)^{m-\ell} \delta_{p,\ell}$ is diagonal and 
$S_m$ is the forward shift on $\mathbb C^{m+1}$ with weight sequence 
$\{1,\ldots ,m\}$, that is, $(\!\!(S_m)\!\!)_{\ell,p}= \ell \delta_{p+1,\ell}$, 
$0\leq p, \ell \leq m$.  Substituting \eqref{Jexp} into \eqref{Ktransrule} and polarizing 
we obtain 
\begin{equation} \label{Kexp}
\bm B^{(\lambda, \m)}(z,\omega) = (1-z \bar{\omega})^{-2\lambda-m} D(z\bar{\omega})
\exp(\bar{\omega}S_m) \bm B^{(\lambda, \m)}(0,0) \exp(z S_m^*) D(z\bar{\omega}). 
\end{equation}

In general, let  $\mathcal{H}$ be a Hilbert space consisting of holomorphic
functions on the open unit disc $\D$ with values in $\C^{m+1}$.  Assume
that $\mathcal{H}$ possesses a reproducing kernel $K:\D\times \D \to
\C^{(m+1)\times (m+1)}$.  The set of vectors 
$\mathcal H_0 = \{K_\omega\xi: \omega\in \bb D, \xi \in \bb C^{m+1}\}$ 
span the Hilbert space
${\mathcal{H}}$.  On the dense set of vectors ${\mathcal H}_0$, we define a map
$T$ by the formula $T K_\omega\xi = \bar{\omega}
K_\omega\xi$ for $\omega\in
\D$.  The following Lemma gives a criterion for boundedness of the
operator $T$.  
\begin{Lemma} \label{bounded}
The densely defined operator $T$ is bounded if and only if for some
positive constant $c$ and for all $n\in \N$
$$ 
\sum_{i,j=1}^n \inner{(c-\omega_j\bar{\omega}_i)
K(\omega_j,\omega_i)x_i}{x_j} \geq 0
$$
for $x_1, \ldots, x_n \in \C^{m+1}$ and $\omega_1, \ldots, \omega_n \in \D$.  If
the map $T: {\mathcal H}_0 \to {\mathcal H}_0
\subseteq {\mathcal{H}}$ is bounded then it is the adjoint of the
multiplication operator on ${\mathcal{H}}$.
\end{Lemma}

The proof is well-known and easy in the scalar case.
We omit the obvious modifications required in the general case.  

It is known and easy to verify that for every $\epsilon > 0$, 
the multiplication operator $M^{(\epsilon)}$, defined by 
$\big (M^{(\epsilon)} f \big )(z) = z f(z)$, is bounded on
$A^{(\epsilon)}$.  Consequently, 
the kernel $B^{\epsilon}$ satisfies the positivity condition of the 
Lemma above for $\epsilon>0$.  Fix $m\in 
\N$.  Consider the reproducing kernel $\bm{B}^{(\lambda,\m)}$.   We recall
that $\bm{B}^{(\lambda, \m)}$ is a positive definite kernel on the unit disc
$\D$ if and only if $\lambda > m/2$.  
\begin{Theorem} \label{}
The multiplication operator $M^{(\lambda, \m)}$ on the Hilbert space 
$\bm A^{(\lambda,\m)}$ is bounded for all $\lambda > m/2$.
\end{Theorem}
\begin{proof}
Let $\epsilon$ be a positive real number such that $\lambda -
\epsilon > m/2$.  Let us find $\m^\prime$ with $\mu_j^\prime > 0$, 
$0\leq j \leq m$, such that 
\begin{equation} \label{muprime}
\bm{B}^{(\lambda, \m)}(z, \omega)  = (1-z\bar{\omega})^{-2\epsilon} \bm{B}^{(\lambda-\epsilon, \m^\prime)}
(z,\omega).  
\end{equation}
Since the multiplication operator is bounded on the Hilbert space
whose reproducing kernel is $(1-z\bar{\omega})^{-2\epsilon}$ for every
$\epsilon > 0$, it follows that we can find $r > 0$ such that 
$(r-z\bar{\omega})(1-z\bar{\omega})^{-2\epsilon}$ is
positive definite. Assuming the existence of $\m^\prime$ as above, 
we conclude that $(r-z\bar{\omega}) \bm{B}^{(\lambda, \m)}(z,\omega)$ 
is positive definite finishing the proof. To find such a $\m^\prime$, it is enough to prove 
$\bm{B}^{(\lambda, \m)}(0,0)= \bm{B}^{(\lambda-\epsilon, \m^\prime)}(0,0)$,  
because then \eqref{B(z,z)} and \eqref{Jexp} (or \eqref{Kexp}) 
immediately imply \eqref{muprime}.

By \eqref{B(0,0)}, writing $L(\lambda)_{\ell j} = \bm B^{(\lambda_j)}(0,0)_{\ell \ell}$, the question becomes whether we can find positive numbers $\mu_j^\prime$ satisfying the equations
\begin{equation} \label{L}
\sum_j L(\lambda)_{\ell j}\, \mu_j^2 = \sum_jL(\lambda-\epsilon)_{\ell j}\, {\mu_j^\prime}^2.
\end{equation}
By \eqref{Bj(0,0)} each $L(\lambda)_{\ell j}$ is continuous in $\lambda$; also 
 $L(\lambda)_{\ell j}=0$ for $\ell < j$, and $L(\lambda)_{0 0} = 1$. It follows that for small 
 $\epsilon > 0$, the system \eqref{L} has solutions satisfying ${\mu_0^\prime}^2=1$, 
 ${\m_j^\prime}^2 > 0$  $(1 \leq j \leq m)$.
\end{proof}

Next we compute the matrix of $M$ with respect to the orthonormal basis 
$\{\mu_j e^j_{n}(z): n \geq 0; 0\leq j \leq m \}$.  
Let $\scr H(n)$ be the linear span of the vectors 
$\{e^j_{n-j}(z): 0 \leq j \leq {\rm min}(m,n) \}$.  
It is clear that $M$ maps the space $\scr H(n)$ into $\scr H(n+1)$. 
(The subspace $\scr H(n)$ of $\bm A^{(\lambda,\m)}(\bb D)$ 
is a ``K-type'' of the \rep $U$.)   We therefore have 
$$
M \mu_j e^j_{n-j} = \sum_{k=0}^m M(n)_{k,j}\, \mu_k\, e^k_{n+1-k}.
$$
Let $E(n)$ be the matrix, determined by \eqref{onb bmA_j}, such that 
$\big( e_{n-j}^j(z)\big)_\ell =E(n)_{\ell,j}z^{n-\ell}$, 
$n \geq j,\: 0\leq j \leq m$.  In this notation, 
$$
E(n)_{\ell,j} \mu_j = \sum_{k=0}^m M(n)_{k,j} E(n+1)_{\ell,k} \mu_k.  
$$
In matrix form, this means 
\begin{align*}
E(n) D(\m) &= E(n+1) D(\m) M(n),\:\mbox{which gives}\\
M(n) &= D(\m)^{-1} E(n+1)^{-1} E(n) D(\m), 
\end{align*}
where $D(\m)$ is the diagonal matrix with $D(\m)_{\ell,\ell}=\mu_\ell$. 
(These are the blocks of $M$ regarded as a ``block shift'' with respect to 
the orthogonal decomposition 
of $\bm A^{(\lambda, \m)}(\bb D) = \oplus_{n=0}^\infty \scr H(n)$.) 

To get information about $M(n)$, we note that, as $n\to \infty$, 
Stirling's formula gives, for any fixed $b\in \bb R$,
$$
\Gamma(n+b) \sim \sqrt{2 \pi} (n+b)^{n+b-1/2} e^{-(n+b)} 
\sim \sqrt{2 \pi} n^{n+b-1/2}(1+\tfrac{b}{n})^n e^{-(n+b)} \sim e^bn^{n+b-1/2}.
$$ 
Applying this we immediately get, by (\ref{onb bmA_j}), 
$$
E(n)_{\ell,j} \sim  n^\ell n^{\lambda - m/2 - 1/2} E_{\ell,j},
$$
where $E$ is the matrix with entries
$$
E_{\ell,j} = \begin{cases} \binom{\ell}{j} \frac{\sqrt{\Gamma(2 \lambda - m + 2 j)}}
{\Gamma(2 \lambda - m + \ell + j)} & \ell \geq j \\
0 & \ell < j \end{cases}
$$
independent of $n$.  Using the diagonal matrix $\bm d(n)$ with 
$\bm d(n)_{\ell,\ell} = n^\ell$, we can write 
$$
E(n)\sim n^{\lambda - m/2 -1/2} \bm d(n) E.
$$ 
It follows that 
\begin{align*}
M(n) &= D(\m)^{-1}E(n+1)^{-1}E(n) D(\m)\\ 
&\sim \big (\frac{n}{n+1}\big )^{\lambda-m/2-1/2} D(\m)^{-1}E^{-1}
\bm d(n+1)^{-1}\bm d(n) E D(\m). 
\end{align*}
Since $\frac{n}{n+1}= 1 + O(\frac{1}{n})$, this implies 
$$
M(n)=I+O(\tfrac{1}{n}),
$$
where $I$ is the identity matrix of order ${m+1}$ and $O(\frac{1}{n})$ 
stands for a $(m+1)\times (m+1)$ matrix each of whose entries is 
$O(\frac{1}{n})$.  

We denote by $U_+$ the operator on $\bm A^{(\lambda,\m)}(\bb D)$ defined by 
$U_+e^j_{n-j} = e^j_{n+1-j}$ ($0 \leq j \leq \mbox{min}(m,n)$, $n-j \geq 0$). 

\begin{Theorem}
The operator $M$ on $\bm A^{(\lambda, \m)}(\bb D)$ is the sum of 
$U_+$ and of an operator in the 
Hilbert-Schmidt class.  In particular, $M$ is bounded and its adjoint 
belongs to the Cowen-Douglas class.
\end{Theorem}
\newsection
\section{Irreducibility}
Let ${\mathcal H}_1$ and ${\mathcal H}_2$ be two \rk Hilbert spaces 
consisting of holomorphic functions on $\D$ taking values in $\C^{m+1}$. 
Suppose that the multiplication \op $M$ on these two Hilbert spaces 
are bounded. Furthermore, assume that the standard set of $m+1$ orthonormal vectors 
$\varepsilon_0,\ldots, \varepsilon_m$ in $\C^{m+1}$, thought of as constant functions on $\D$, are in 
both $\mathcal H_1$ and $\mathcal H_2$. Since 
$
\big ( \sum_{i=0}^m p_i(M) \varepsilon_i \big )(z) = \sum_{i=0}^m p_i(z) \varepsilon_i 
$
for polynomials $p_i$ with scalar coefficients, it follows  that the polynomials 
$\bm p(z) = \sum_{i=0}^m p_i(z) \varepsilon_i$
belong to these Hilbert spaces.  We assume that the polynomials $\bm p$ are dense
in both of these Hilbert spaces.

Suppose that there is a bounded operator $X:\scr H_1 \to \scr H_2$ satisfying $M X = X M$. Then    
$$
(X\bm p)(z) = \big ( X \sum_{i=0}^m p_i  \varepsilon_i \big )(z)
= \big ( X \sum_{i=0}^m p_i(M) \varepsilon_i \big )(z)
=\big (\sum_{i=0}^m p_i(M) X \varepsilon_i \big )(z)
=\big (\sum_{i=0}^m p_i(M) X \varepsilon_i \big )(z).
$$
Now, if we let $(X\varepsilon_i)(z) = \sum_{j=0}^m x^j_i(z) \varepsilon_j$,   
then  $(X\bm p)(z) = \Phi_X(z) \bm p(z)$, where $\Phi_X(z)= (\!(x^j_i(z))\!)_{j,i=0}^m$.  Since the 
polynomials $\bm p$ are dense, it follows that  $(Xf)(z) = \Phi_X(z) f(z)$ for all $f \in \mathcal H_1$. 

We calculate the adjoint of the intertwining operator $X$.  We have 
\begin{align*}
\inner{XK_1(\cdot,\omega)\xi}{K_2(\cdot,u) \eta} &= 
\inner{\Phi_X(\cdot)K_1(\cdot,\omega)\xi}{K_2(\cdot,u) \eta} = \inner{\Phi_X(u)K_1(u,\omega)\xi}{\eta}\\
&=\inner{K_1(u,\omega)\xi}{\overline{\Phi_X(u)}^{\rm tr} \eta} =
\inner{K_1(\cdot,\omega)\xi}{K_1(\cdot,u) \overline{\Phi_X(u)}^{\rm tr} \eta}
\end{align*}
for all $\xi,\eta \in \C^{m+1}$, that is, 
\begin{equation} \label{Xadj}
X^* K_2(\cdot,u) \eta = K_1(\cdot,u) \overline{\Phi_X(u)}^{\rm tr} \eta,
\end{equation}
for all $\eta \in \C^{m+1}$ and $u\in \D$. Hence the intertwining operator $X$ is unitary if and only if   
there exists an invertible holomorphic function $\Phi_X:\D \to \bb C^{(m+1)\times (m+1)}$  
satisfying 
\begin{equation} \label{unitaryrel}
K_2(z,\omega) = \Phi_X(z) K_1(z,\omega) \overline{\Phi_X(\omega)}^{\rm tr}. 
\end{equation}

Let $\mathcal{H}$ be a \hs consisting of $\bb{C}^n$ - valued holomorphic functions on 
$\D$. Assume that $\mathcal{H}$ has a reproducing kernel, say $K$.  Let  $\Phi$  be a $n\times n$
invertible matrix valued holomorphic function on $\D$ which is invertible.
For $f\in \mathcal{H}$, consider the map $X: f\mapsto
\tilde{f}$, where
$\tilde{f}(z)= \Phi(z)f(z)$.  Let $\tilde{\mathcal{H}}=\{\tilde{f}:
f\in \mathcal{H}\}$.  The requirement that the map $X$ is
unitary, prescribes a \hs structure for the function space
$\tilde{\mathcal{H}}$.  The \rk for $\tilde{\mathcal{H}}$ is clearly 
\begin{equation} \label{rk for tilde}
\tilde{K}(z,\omega) = \Phi(z) K(z,w) \Phi(\omega)^*.
\end{equation}
It is easy to verify that $X M X^*$ is
the multiplication \op $M: \tilde{f} \mapsto z \tilde{f}$ on the \hs
$\tilde{\mathcal{H}}$.  Suppose we have a unitary \rep $U$ given by a
multiplier $J$ acting on $\mathcal{H}$ according to (\ref{Ktransrule}).   
Transplanting this action to $\tilde{\mathcal{H}}$ under
the isometry $X$, it becomes
$$
\big (\tilde{U}_{g^{-1}} \tilde{f}\big )(z) = \tilde{J}(g, z)\tilde{f}
(g\cdot z),
$$
where the new multiplier $\tilde{J}$ is given in terms of the original
multiplier $J$ by
\begin{equation} \label{eqmult}
\tilde{J}(g,z) = \Phi(z)J(g,z) \Phi(g\cdot z)^{-1}.
\end{equation}
Of course, now $\tilde{K}$ transforms according to (\ref{Ktransrule}),
with the aid of $\tilde{J}$.

\begin{Lemma} \label{redcond}
Suppose that the \op $M$ acting on the \hs $\mathcal H$ with \rk $K$ is bounded, 
the constant vectors $\varepsilon_0, \ldots, \varepsilon_m$ are in $\mathcal H$, 
and that the polynomials $\bm p$ are dense in $\mathcal H$.
If there exists a (self adjoint) projection $X$ commuting with the \op
$M$ then  
$$
\Phi_X(z) K(z,\omega) = K(z,\omega)\overline{\Phi_X(\omega)}^{\rm tr}    
$$
for some holomorphic function $\Phi_X: \bb D \to \bb C^{(m+1)\times (m+1)}$ with $\Phi_X^2 = \Phi_X$.  
\end{Lemma}
\begin{proof}
We have already seen that any such \op $X$ is multiplication by a 
holomorphic function $\Phi_X$. To complete the proof, note that 
$$ 
\Phi_X(\cdot) K(\cdot,\omega)\xi = X K(\cdot,\omega) \xi = X^*K(\cdot,\omega) \xi=  
K(\cdot,\omega)\overline{\Phi_X(\omega)}^{\rm tr} \xi
$$
for all $\xi \in \C^{m+1}$.
\end{proof}

From the Lemma,  putting $\omega=0$, we see that $\Phi_X(z) = K(z,0)\overline{\Phi(0)}^{\rm
tr}K(z,0)^{-1}$ for any self adjoint intertwining operator $X$.  Furthermore, $X_0:=\Phi_X(0)$ 
is an ordinary projection on $\bb C^{m+1}$, if $K(0,0)=I$.  The multiplication operator on the two Hilbert
spaces $\scr H$ with \rk $K$ and $\scr H_0$ with \rk $K_0(z,\omega)=
K(0,0)^{-1/2}K(z,\omega)K(0,0)^{-1/2}$ are unitarily equivalent via
the unitary map $f \mapsto K(0,0)^{-1/2} f$. The \rk $K_0$ has the
additional property that $K_0(0,0)=I$.  Therefore, we conclude
that $M$ is reducible if and only if there exists a projection $X_0$ on
$\bb C^{m+1}$ satisfying
\begin{equation} \label{Mirr}
X_0 K_0(z,0)^{-1}K_0(z,\omega) K_0(0,\omega)^{-1} = 
K_0(z,0)^{-1} K_0(z,\omega) K_0(0,\omega)^{-1} X_0. 
\end{equation}
This is the same as requiring the existence of a projection $X_0$ which 
commutes with all the coefficients in the power series expansion of the 
function  
$\widehat{K_0}(z,\omega):= K_0(z,0)^{-1}K_0(z,\omega) K_0(0,\omega)^{-1}$ 
around $0$. We also point out that $\widehat{K_0}$ is the normalized kernel 
in the sense of \cite{C-S} and is characterized by the property 
$\widehat{K_0}(z,0)\equiv 1$.

For the rest of this section, we set $B:= \bm B^{(\lambda,\m)}(0,0)$ and 
$S:=S_m$, as in Section $4$.
\begin{Lemma} \label{MirrB}
The operator $M:=M^{(\lambda,\m)}$ on the \hs $A^{(\lambda, \m)}$ is \irr 
if and only if there is no projection $X_0$ on $C^{m+1}$ commuting with 
all the coefficients in the power series expansion of the function 
$$
(1-z\bar{\omega})^{-2\lambda-m} B^{1/2}
\exp{(-zS^*)}B^{-1} D(z\bar{\omega})\exp{(\bar{\omega}S)}B
\exp(zS^*)D(z\bar{\omega})B^{-1}\exp(-\bar{\omega}S)B^{1/2}, 
$$ 
around $0$.  
\end{Lemma}

\begin{proof}
From \eqref{Kexp}, we have  
$
\bm B^{(\lambda, \m)}_0(z,0) = B^{1/2}\exp(z S^*)B^{-1/2}, 
$
where $\bm B_0^{(\lambda,\m)}:=B^{-1/2} \bm B^{(\lambda,\m)}B^{-1/2}$.  
To complete the proof, using \eqref{Kexp}, we merely verify that 
\begin{eqnarray*}
\lefteqn{\widehat{\bm B_0}(z,\omega)}\\ 
&=&\big (\bm B_0^{(\lambda,\m)}(z,0) \big )^{-1} 
\bm B_0^{(\lambda,\m)}(z,\omega)\big (\bm B_0^{(\lambda,\m)}
(0,\omega) \big )^{-1}\\
&=& (1-z\bar{\omega})^{-2\lambda-m} B^{1/2}
\exp{(-zS^*)}B^{-1} D(z\bar{\omega})\exp{(\bar{\omega}S)}B
\exp(zS^*)D(z\bar{\omega})B^{-1}\exp(-\bar{\omega}S)B^{1/2}.
\end{eqnarray*}
\end{proof}

Let $D_s$ denote the coefficient of $(-1)^sz^s\bar{\omega}^s$ in the 
expansion of $D(z\bar{\omega})$ and $\tilde{D}_s = B^{-1} D_s$. 
(The choice of $D_s$ ensures that the diagonal sequence in 
$\tilde{D}_s$ is positive.)     

\begin{Lemma} \label{product}
If $(\!\!( {S^*}^i \tilde{D}_s S^p B {S^*}^q \tilde{D}_t S^j)\!\!)_{k n} 
\not = 0$ for some choice of $i,j,s,t,p,q$ in $\{0,1,\ldots ,m\}$ then 
\begin{eqnarray*}
0\leq s \leq m-k-i, & 0\leq t \leq m-n-j;\\
0\leq p \leq k +i, & 0 \leq q \leq n+j;
\end{eqnarray*}
and $k+i-p = n+j -q$. 
\end{Lemma}
\begin{proof} By the definition of $S$ we have 
$$ S^p:\begin{cases} e_\ell \mapsto (l+1)_p e_{\ell+p} & \mbox{\rm ~if~} 
0\leq \ell \leq m-p\\
e_\ell \mapsto 0 & \ell > m-p, 
\end{cases}
$$
Also, 
$$\tilde{D}_s:\begin{cases} e_\ell \mapsto c\, e_\ell & 
\mbox{\rm ~if~}0\leq \ell \leq m-s\\
e_\ell \mapsto 0 & \ell \geq m-s+1,
\end{cases}
$$
where $c$ is a non-zero constant depending on $\ell, s$. Therefore 
$$Q:={S^*}^i\tilde{D}_s S^p:\begin{cases} e_\ell \mapsto c^\prime e_{\ell+p-i} & 
\mbox{\rm ~if~} 0\leq i \leq m-p-s \mbox{\rm ~and~} \ell+p-i \geq 0\\
e_\ell \mapsto 0 & \mbox{\rm otherwise} 
\end{cases}
$$
for some non-zero constant $c^\prime$. Hence the full condition 
for $Q_{k \ell}\not = 0$ is 
\begin{equation} \label{cond1}
i-p \leq \ell \leq m -p-s,\, k = \ell + p - i.
\end{equation}
Let $R:={S^*}^q\tilde{D}_tS^j$.  By what we have just proved, it
follows that $R_{\ell n}\not = 0$ if and only if 
\begin{equation} \label{cond2}
q-j\leq n \leq m-j-t,\, \ell= n + j - q.
\end{equation}
Now since $B$ is diagonal non-zero, we have $(QBR)_{kn} 
\not = 0$ if and only if $Q_{k \ell},\,R_{\ell n} \not = 0$ for some $\ell$.  
By  (\ref{cond1}) and (\ref{cond2}) this happens exactly when  
\begin{equation} \label{condf}
0 \leq \ell = k+i-p = n+j-q = \ell \leq m,\, k+i \leq m -s \mbox{\rm ~and~} n+j
\leq m-t.
\end{equation}
\end{proof}
Let $a(\ell)$ denote the coefficient of 
$z^{m+\ell+1}\bar{\omega}^{m+\ell}$ ($0 \leq \ell < m-1$) in the polynomial $A$, where  
\begin{eqnarray*}
A(z,\omega)&=& \exp{(-zS^*)}B^{-1} D(z\bar{\omega})\exp{(\bar{\omega}S)}B
\exp(zS^*)D(z\bar{\omega})B^{-1}\exp(-\bar{\omega}S)\\
&=& \sum (-1)^i\frac{{S^*}^i}{i!}z^i (-1)^s\tilde{D}_s z^s\bar{\omega}^s 
\frac{S^p}{p!}\bar{\omega}^p B \frac{{S^*}^q}{q!}z^q(-1)^t\tilde{D}_t 
z^t\bar{\omega}^t (-1)^j\frac{S^j}{j!}\bar{\omega}^j,   
\end{eqnarray*}
where the sum is over $0 \leq i,j,p,q,s,t \leq m$.
\begin{Lemma} \label{zwproduct} 
For $0\leq \ell \leq m-1$, 
$a(\ell)_{k,n}= \begin{cases} \mbox{\rm not zero} & 
\mbox{\rm ~if~} k = m-\ell-1 \mbox{\rm ~and~}n= m-\ell\\ 
\mbox{\rm zero} & \mbox{\rm ~if~} k-n \not = 1 \mbox{\rm ~or~} 
k > m -\ell - 1\end{cases}$. 
\end{Lemma}
\begin{proof} Clearly, $A(z,w) = \sum A_{ijpqst}
  z^{i+s+q+t}\bar{\omega}^{s+p+t+j}$, where the sum is over $0\leq
  i,j,p,q,s,t \leq m$.  Therefore, $a(\ell) = \sum c {S^*}^i
  \tilde{D}_s S^p B {S^*}^q\tilde{D}_t S^j$, where the sum is over all
  $i,j,p,q,s,t$ such that $s+t+i+q = m+\ell+1$ and $s+t+p+j=m+\ell$;
  $c=\frac{(-1)^{i+j+s+t}}{i!\,j!\,p!\,q!}$.

  It follows from the preceding Lemma that if $a(\ell)_{k,n} \not =
  0$, then $i-j+q-p = n-k$.  However, for the terms occuring in the
  sum, we now have $i-j+q-p = (s+t+i+q) - (s+t+p+j)=1$. Thus if
  $a(\ell)_{k,n} \not = 0$ then $n-k=1$.

Furthermore, if $a(\ell)_{k,n} \not = 0$, then we also have $m+\ell+1
= (s+t+i+q)$. Hence $m+\ell+1 - (s+t+i) = q \leq n+j$ from the last
inequality of the preceding Lemma, that is, $s+t+i+j \geq m+\ell+1
-n$.  This along with $s+t+i+j \leq 2m -k -n$, which is obtained by
adding the first two inequalities of the preceding Lemma, gives $k\leq
m-\ell-1$

The proof of the second part of the Lemma is now complete.  

If $k=m-\ell-1$ and $n=m-\ell$, for the terms occuring in the sum for
$a(\ell)$, we have $s+t+i+j=2\ell+1$.  It follows that
$a(\ell)_{m-\ell-1,m-\ell}$ is a sum of negative numbers.  This proves
the first part of the Lemma.
\end{proof}

\begin{Theorem} \label{}
The multiplication \op $M:=M^{\lambda,\m}$ on the \hs $\bm A^{(\lambda,\m)}$ 
is irreducible.
\end{Theorem}

\begin{proof}
Suppose there exists a non-trivial projection 
$P$ commuting with $\widehat{\bm B_0}(z,\omega)$ for all $z,\omega\in \bb D$. 
Then, by Lemma \ref{MirrB} $P$ must commute with $B^{1/2}A(z,\omega)B^{1/2}$
for all $z,\omega \in \bb D$.  By Lemma \ref{zwproduct}, $P$ must commute with 
all shifts, that is, matrices $T$ such that $T_{k n} = 0$ unless $k-n=1$.   
This is a contradiction.
\end{proof}
\newsection
\section{Inequivalence}
Let ${\rm pr} :E_T \to \D$ be the holomorphic vector bundle 
corresponding to an operator $T\in {\rm B}_k(\D)$.  The operator $T$ is 
homogeneous if and only if for any $g\in G$, there exists an automorphism 
$\hat{g}$ of the bundle $E_T$ covering $g$, that is, the diagram 
$\begin{CD}
E_T @>\hat{g}>> E_T\\@VV{\rm pr}V @VV{\rm pr}V\\\D @>g>>\D
\end{CD}\:\:$ 
is commutative.  

\begin{Theorem} \label{tildeGacts}
If $T$ is a \hom \op in ${\rm B}_k(\D)$ then the 
the universal covering group $\tilde{G}$ of $G$ acts 
on $E_T$ by automorphisms.
\end{Theorem}
\begin{proof}
Let $\hat{G}$ be the group of automorphisms of $E_T$.  This is a Lie
group.  Let $p:\hat{G} \to G$ be the natural homomorphism.  Let
$N=\ker p$, the automorphisms fixing all the points of $\D$.  Then
$\hat{G}/N\simeq G$, and for the corresponding Lie algebras, we have
$\hat{\ger{g}}/\ger{n} \simeq\ger{g}$.  Since $\ger g$ is
semisimple, by the Levi decomposition, there is a subalgebra
$\hat{\ger g}_0 \subseteq \hat{\ger g}$ such that $\hat{\ger g} =
\hat{\ger g}_0 +\ger n$, where the sum is a vector space direct   
sum.  Let $\hat{G}_0$ be the corresponding analytic subgroup.  

There is a neigbourhood $U$ of $e \in \hat{G}_0$ such that $p_{| U}$
is a homeomorphism onto a neighbourhood $p(U)$ of $e\in G$.  But then
$p(\hat{g}U) = p(\hat{g})p(U)$.  So, $p$ is a homeomorphism of a
neighbourhood of any point $\hat{g}\in \hat{G}_0$ to a neighbourhood
of $p(\hat{g})$ in $G$. It follows that the image of $p$ is an open
subgroup and so must equal $G$. Therefore, $\hat{G}_0$ is a covering
group of $G$.

Now, $\hat{G}_0$ acts on $E_T$ by automorphisms and projects to $G$.
The universal cover $\tilde{G}$ now also acts on $E_T$.

\end{proof}
\begin{Remark}
With slightly more work one can prove that $\hat{\ger g}_0$ is an 
ideal and therefore the $\tilde{G}$ action on $E_T$ is unique.  
If $T$ is \irr it is known independently (cf. \cite{surv}) that the 
$\tilde{G}$ action is unique.
\end{Remark}

\begin{Theorem} \label{}
For every $m \geq 1$, the operators $M^{(\lambda,
\m)},\, \lambda > \frac{m}{2}; \mu_1, \ldots, \mu_m >0$ 
are mutually unitarily inequivalent.
\end{Theorem}
\begin{proof}
Suppose $M^{(\lambda,\m)}$ and $M^{(\lambda^\prime,\m^\prime)}$ are unitarily equivalent. 
Then the corresponding Hermitian holomorphic bundles are isomorphic \cite{C-D}.  Now, 
by \eqref{unitaryrel} and \eqref{eqmult}, multipliers $J$ and $J^\prime$ giving the 
$\tilde{G}$ action on $\mathbb A^{(\lambda,\m)}$ and $\mathbb A^{(\lambda^\prime,\m^\prime)}$ 
are equivalent in the sense that there exists a invertible matrix function $\phi(z)$, 
holomorphic in $z$, such that  
$$
\Phi(z) J(g,z) \Phi(g z)^{-1} = J^\prime(g,z)
$$
on $\tilde{G} \times \D$. 
Setting here $g=p_{-z}$, \eqref{Jexp} gives 
$$
\Phi(z)= (1-|z|^2)^{\lambda - \lambda^\prime} D(|z|^2) \exp(-\bar{z}S_m)D(|z|^2) F(0) \exp(\bar{z}S_m)D(|z|^2)^{-1}.
$$
The right hand side is real anlytic in $z,\bar{z}$ on $\D$.  Since $\Phi$ is holomorphic, $\Phi(z) = \Phi(0)$ 
identically.  Looking at the Taylor expansion, we obtain
$$
S_m\Phi(0) = \Phi(0) S_m. 
$$
This implies that $\Phi(0) = p(S_m)$, a polynomial in $S_m$. 
(Note that $S_m$ is conjugate to $S$, the unweighted shift with entries 
$S_{\ell p}= \delta_{p+1\,\ell}$, which is its Jordan canonical form.   
For $S$ the corresponding property is easy to see.)   
We write 
$$
D^{1,1} = \left . \frac{\partial^2}{\partial z \partial \bar{z} } \right |_{0}\, D(|z|^2) = - 
\left ( \begin{smallmatrix} m & & & & \\  & m-1 &  &  \\  &  &\ddots & & \\
 & & & 1 & \\  &  &  & & 0\end{smallmatrix} \right ),
$$
and for the Taylor coefficient of $z \bar{z} = |z|^2$ we obtain 
$$
(\lambda -\lambda^\prime) \Phi(0) + D^{1,1} \Phi(0) - \Phi(0) D^{1,1} = 0.
$$
Consider the diagonal of this matrix equality. All diagonal elements of 
$\Phi(0) = p(S_m)$ are the same number $x\not = 0$ (since $p(S_m)$ is triangular and invertible).  
Hence $\lambda - \lambda^\prime = 0$.  Now, since the diagonal entries of $D^{1,1}$ are all different, 
$\Phi(0)$ must be diagonal.  So, $\Phi(0)= x I_{m+1}$.  Also, 
$\Phi(0)$ intertwines the operators  $M^{(\lambda,\m)}$ and $M^{(\lambda^\prime,\m^\prime)}$, hence 
$\Phi(0) \bm B^{(\lambda,\m)}(z, \omega) \Phi(0)^* = \bm B^{(\lambda^\prime,\m^\prime)}(z,\omega)$ 
as in \eqref{unitaryrel}.  Using this with $z=\omega=0$ and using \eqref{Bj(0,0)}, \eqref{B(0,0)} 
we get $|x|^2 \mu_j^2 = {\mu_j^\prime}^2$ for all  $j$.     Since $\mu_0=1=\mu_0^\prime$, it follows 
that $|x|^2 = 1$ and $\mu_j = \mu_j^\prime$ for $1\leq j \leq m$.   
\end{proof}

\section{Some further remarks}
We presented the \ops $M^{(\lambda,\m)}$ in as elementary a way
as possible, but this presentation hides the natural way in which
these \ops can be found to begin with.  One such way, which was
actually followed by the authors, is to start with an irreducible
finite dimensional \rep $\varrho_m$ of ${\rm SU}(1,1)$ (it is well
known that there is exactly one for every natural number $m$), observe
that $J(g,z)=\varrho_m(g^{-1})$ can be used as a multiplier, then
transform this multiplier to a more convenient form, to construct a
\rep of $\tilde{G}$ and proceed from there.  This was also the
procedure of Wilkins \cite{W} who worked with the identical ($2$-
dimensional) \rep of ${\rm SU}(1,1)$. The authors are planning to
write an expository article in which there will be some details of
this approach.

Another way, which is probably the most natural one, is to start with
the process of holomorphic induction to construct the \hom vector
bundles which are to be the \chd bundles of our \Ops. It is well known
that every finite dimensional \rep $\varrho$ of the triangular
subalgebra $\ger t$ of $\mathfrak{sl}(2, \bb C)$ (the Lie algebra of the
stabilizer of $0$ in ${\rm SL}(2, \bb C)$ acting on the extended
complex plane) gives rise to a $\tilde{G}$-\hom holomorphic vector
bundle, from which the $\varrho$ can be reconstructed.  Refining this
statement, it is easy to see that the $\tilde{G}$-\hom holomorphic
Hermitian vector bundles are in one-to-one correspondence with the
unitary equivalence classes of \reps $\varrho$ of $\ger t$ on on the
finite dimensional Hilbert spaces $\bb C^n$ with the added property
that $\varrho$ is skew Hermitian on the (one-dimensional;) subalgebra
$\ger k$, the Lie algebra of the stabilizer of $0$ in ${\rm SU}(1,1)$.

If we start with the restriction to $\ger t$ of the
$(m+1)$-dimensional \rep $\varrho_m$ of $\mathfrak{sl}(2, \C)$ and we put
on the \rep space all possible inner products so that the
requirement concerning $\ger k$ is satisfied then we obtain a family
of bundles parametrized by $\lambda \in \R$ and $\mu_1, \ldots , \mu_m
>0$.  It can then be shown that these bundles correspond to \chd \ops
if and only if $\lambda > m/2$, and in this case the corresponding \op is
$M^{(\lambda,\m)}$.

One can use this approach starting with any finite dimensional \rep
$\varrho$ of $\ger t$.  Such a $\varrho$ can always be written as
$\varepsilon_\lambda \otimes \varrho_0$, where $\varepsilon_\lambda$
($\lambda \in \R$) is a one dimensional \rep of $\ger k \cong \R$
extended trivially to $\ger t$ and $\varrho_0$ is normalized in a
certain way.  There is always a corresponding \hom Hermitian vector
bundle and a number $\lambda_\varrho$ such that for $\lambda >
\lambda_\varrho$ the bundle corresponds to a \hom \chd \Op.

In this generality one cannot expect as explicit results as in the
present paper, but one can proceed to still make fairly precise
statements.  In this way one gets a kind of classification of all \hom
\chd \Ops. This will be the subject of a second article in this series.  

Finally we mention that many of our arguments extend without change to
the case of operator tuples and holomorphic vector bundles over
bounded symmetric domains in several complex variables. There are, of
course, a number of new features (cf. \cite{twist,Arazy-Z}) as well in this 
general situation which still have to be explored in the future.

\bibliographystyle{amsplain}
\providecommand{\bysame}{\leavevmode\hbox to3em{\hrulefill}\thinspace}
\providecommand{\MR}{\relax\ifhmode\unskip\space\fi MR }
\providecommand{\MRhref}[2]{%
  \href{http://www.ams.org/mathscinet-getitem?mr=#1}{#2}
}
\providecommand{\href}[2]{#2}

\end{document}